# Mathematical Modeling of Carbon Dioxide Emissions with GDP Linkage: Sensitivity Analysis and Optimal Control Strategy

**Hua Liu[1*], Zhuoma Gangji[1], Yumei Wei[2], Jianhua Ye[3], Gang Ma[1]**
1. School of Mathematics and Computer Science, Northwest Minzu University, Lanzhou Gansu 730000, China
2. Experimental Teaching Deparment, Northwest Minzu University, Lanzhou, Gansu 730030, China
3. School of Preparatory Education, Northwest Minzu University, Lanzhou, Gansu 730000, China
* Corresponding author: jslh@xbmu.edu.cn

**Abstract**

Climate change and global warming are among the most significant issues that humanity is currently facing, and also among the issues that pose the greatest threats to all mankind. These issues are primarily driven by abnormal increases in greenhouse gas concentrations. Mathematical modeling serves as a powerful approach to analyze the dynamic patterns of atmospheric carbon dioxide. In this paper, we established a mathmetical model with four state variables to investigate the dynamic behavior of the interaction between atmospheric carbon dioxide, GDP, forest area and human population. Relevant theories were employed to analyze the system's boundedness and the stability of equilibrium points. The parameter values were estimated with the help of the actual data in China and numerical fitting was carried out to verify the results of the theoretical analysis. The sensitivity analysis of the compartments with respect to the model parameters was analyzed by using the Partial Rank Correlation Coefficient (PRCC) and the Latin Hypercube Sampling test. Apply the optimal control theory to regulate the atmospheric carbon dioxide level and provide the corresponding numerical fitting. Finally, corresponding discussions and suggestions were put forward with the help of the results of the theoretical analysis and numerical fitting.

# 1 Introduction

Since the commencement of the Industrial Revolution, the accelerated pace of industrialization in human society has notably aggravated the emission of greenhouse gases, with carbon dioxide ($CO_2$) emerging as the principal contributor to this phenomenon. This acceleration, driven by exponential growth in fossil fuel combustion, industrial processes, and deforestation, has disrupted the natural carbon cycle, thereby exacerbating the greenhouse effect and altering global climatic patterns (Prabodhi et al. 2020). Throughout this epoch, anthropogenic activities, principally fossil fuel combustion, deforestation, and land-use transformation have propelled atmospheric $CO_2$ concentrations from approximately 280 parts per million (ppm) to over 400 ppm. This unprecedented rise, equivalent to a 43% increase since the pre-industrial era, reflects a systematic disruption of the global carbon cycle, with profound implications for climate stability and ecological balance (ESRL et al. 2018). This increase in concentrations is the main driver of global warming, which exacerbates the rise in global temperatures and triggers widespread climate change phenomena (Rahmstorf et al. 2009). By February 2025, the atmospheric carbon dioxide concentration had surged to an astonishing 426.13 parts per million (ppm) (Lan et al. 2025). Fig. 1 depicts the temporal dynamics of annual global greenhouse gas emissions from 1980 to 2023, revealing a consistent upward trajectory. While inter-annual fluctuations are evident, the secular growth trend remains statistically robust, underscoring the escalating severity of greenhouse gas emissions and the urgent imperatives for climate change mitigation. This trend not only reflects the cumulative impact of anthropogenic activities but also highlights the critical need for concerted, systemic action to address this planetary challenge (Ritchie et al. 2020).

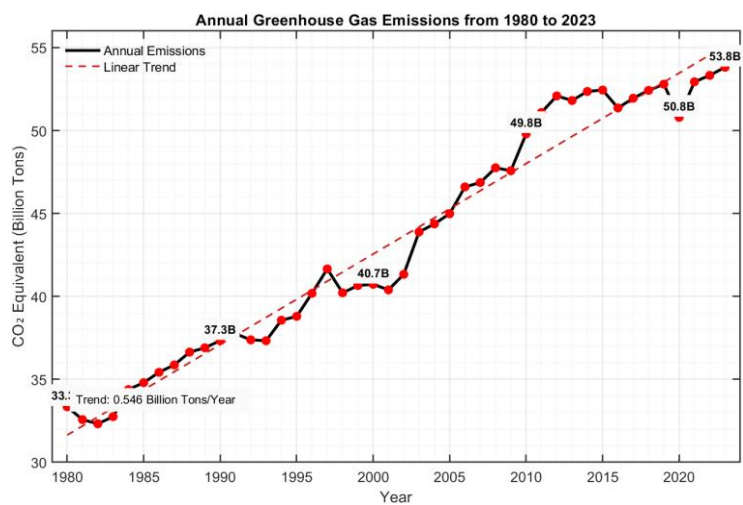

**Fig. 1** Trends of $CO_2$

Climate change is manifested in rising global average temperatures, frequent extreme weather events, accelerating glacier melting, continuous sea-level rise, and significant ecosystem disruptions (Rummukainen 2013). These changes pose a serious threat to global agricultural production, water allocation, biodiversity conservation and human health (Hales et al. 2006). For example, sea-level rise endangers coastal communities and infrastructure, extreme weather intensifies natural disaster risks, and ecosystem disruptions alter species distributions and survival rates (FAO 2020).

In response to global warming and climate change, the international community has taken various actions. The Paris Agreement stands as a key milestone, aiming to keep global temperature rise well below 2°C and pursue efforts for 1.5°C (Blanford et al. 2014). Countries are adopting multiple strategies to cut greenhouse gas emissions: developing renewable energy, boosting energy efficiency, implementing emissions trading schemes, and restoring forests (Çeliktaş et al. 2021). The alterations in extreme weather conditions, coupled with the rise in global surface temperature and climate change, have given rise to the spread of a host of diseases (Prabodhi et al. 2020). Climate change is evidenced by multiple phenomena: an increase in global mean temperature, a surge in extreme weather occurrences, glacial and ice sheet melt, sea-level ascent, and significant ecosystem alterations (Agarwal et al. 2017). Mathematical modeling can effectively visualize the dynamic behavior of atmospheric carbon dioxide, enabling better adoption of corresponding measures to alleviate the level of atmospheric carbon dioxide (Devi 2020).

Numerous mathematical frameworks have been put forward to investigate the impacts of factors such as carbon capture technology, population pressure, reforestation, vehicle $CO_2$ emissions, technology choice and urbanization on the dynamics of atmospheric carbon dioxide (Achimugwu et al. 2023; Jha et al. 2021; Misra et al. 2015; Arora et al. 2024; Jha et al. 2024; Misra et al. 2015; Arora et al. 2024; Jha et al. 2024; Donald et al. 2024; Misra et al. 2018; Misra et al. 2017; Bansal et al. 2024). In particular, Devi and Gupta . (2018) proposed a nonlinear mathematical model to simulate changes in the ability of plants to absorb atmospheric carbon dioxide. The paper indicates that afforestation represents a scientific approach to lower atmospheric carbon dioxide concentrations. In (Misra and Jha. 2022), a three-dimensional mathematical model has been constructed to analyze the impact of budget allocation on the reduction of atmospheric carbon dioxide concentrations. As the level of carbon dioxide in the atmosphere increases, the growth rate of budget allocation may lead to a stability switch through the hopf-bifurcation. Mishra et al. (2019) discussed the use of green belt planting and seaweed farming to reduce atmospheric carbon dioxide ($CO_2$), and model analysis showed that the use of plants for photosynthesis by planting leafy trees in the green belt around the emission source, and through seaweed farming, could effectively reduce atmospheric $CO_2$ levels. Tandon (2023) employed a mathematical model to investigate the impacts of mining activities on the dynamic natural interactions between plants and carbon dioxide. The study revealed that mining activities notably elevated atmospheric carbon dioxide concentrations and caused damage to plants, thereby hindering the system's ability to attain a stable state. Misra and Verma (2023) studied the effects of population and forest biomass on atmospheric carbon dioxide. The results show that when human deforestation exceeds a certain threshold, the system will occur hopf-bifurcation and become unstable. Most of the above-mentioned literatures have considered the kinetic relationship between carbon dioxide and forest. Caetano et al. (2008) focused their analysis on how GDP affects atmospheric carbon dioxide, neglecting to consider the interrelationship between carbon dioxide and forest ecosystems. Our modeling framework accounts for the interactive relationship between $CO_2$ emissions and GDP, constructing a four-dimensional mathematical model to capture the dynamic behavior of atmospheric carbon dioxide.

China is the country with the largest carbon emissions in the world, accounting for nearly one-third of the global total carbon emissions. And this is inextricably linked to China's rapidly growing economy. Therefore, researching the dynamic relationship of the interaction between China's GDP index and the carbon cycle is of far-reaching significance in the efforts to mitigate global warming and climate change.

## 2 Mathematical Model

In this work, we proposed a mathematical model to research the dynamics of $CO_2$ emissions, GDP, forest area and human population. The variables are defined as follows:

(1) $C(t)$: the concentration of atmospheric $CO_2$ (in ppm).

(2) $G(t)$: the gross domestic product (in billion USD$).

(3) $F(t)$: the forest area (in million hectares.).

(4) $N(t)$: the human population (in million).

The atmospheric carbon dioxide emissions stem from two categories: The emission of natural factors (such as volcanic eruptions, respiration processes of plants and animals, etc.), this constant growth term we denote as $\alpha$. The emission of carbon dioxide caused by human activities which is proportional to the population (Newell and Marcus. 1987), we denote this increase as parameter $\phi$. Rapid economic growth is often accompanied by a lot of industrialization (Poterba .1993), i.e. more carbon dioxide is emitted into the atmosphere, which we record as $\beta$. Forest area absorbs atmospheric carbon dioxide through the process of photosynthesis and leads a decrease in atmospheric carbon dioxide (Panja. 2020), which we denote as $\eta$. The lifetime of atmospheric carbon dioxide is usually 30 to 95 years (Jacobson. 2005), we denote the natural loss coefficient of atmospheric carbon dioxide as $p$. Let $\mu$ represent the growth rate of GDP. Economic growth, in turn, can reduce carbon dioxide in the atmosphere through activities such as clean technologies, we use $\varepsilon$ to represent it. Based on these assumptions, the dynamics of atmospheric carbon dioxide are governed by the following equation:

$$\frac{dC}{dt} = \alpha + \phi N + (\beta - \varepsilon)G - \eta CF - pC, \tag{1}$$

$$\frac{dG}{dt} = \mu - \varepsilon G. \tag{2}$$

In the modeling process, we assume that $\omega$ represent the natural growth rate of forest areas and $K$ represent the carrying capacity of forest area. Human population growth often leads to an increase in demand for deforestation, which may include agricultural land, urban sprawl and infrastructure development. To meet these needs, forests may be cut down and leading to a decrease in forest area (Jyotsna. 2024), this anthropogenic deforestation coefficient is denoted as $\theta$. The absorption of the right amount of carbon dioxide will promote the growth of forest area more densely (Gautam. 2024). The growth rate of forest area caused by the absorption of carbon dioxide we denoted as $\sigma$. Based on these assumptions, the dynamics of forest areaare governed by the following equation:

$$\frac{dF}{dt} = \omega F\left(1 - \frac{F}{K}\right) - \theta NF + \eta\sigma CF. \tag{3}$$

We assume that the population follows logistic growth, with $s$ and $M$ denoting the natural growth rate of the population and the carrying capacity of the population respectively. Let $v$ represents the contribution of forest area to the population (e.g. provision of food and resources, conditions, climate, etc.) (Jha .2024). The absorption of carbon dioxide by the human body can have dire consequences, it can directly lead to death or exacerbate the spread of specific diseases, we use parameter $\pi$ to express the rate of loss of the population caused by carbon dioxide (Arora. 2025). Based on these assumptions, the dynamics of popultion governed by the following equation:

$$\frac{dN}{dt} = sN\left(1 - \frac{N}{M}\right) + \theta v NF - \pi CN. \tag{4}$$

In summary, our mathematical model is as follows:

$$\begin{aligned}
\frac{dC}{dt} &= \alpha + \phi N + (\beta - \varepsilon)G - \eta CF - pC, \\
\frac{dG}{dt} &= \mu - \varepsilon G, \\
\frac{dF}{dt} &= \omega F\left(1 - \frac{F}{K}\right) - \theta NF + \eta \sigma CF, \\
\frac{dN}{dt} &= sN\left(1 - \frac{N}{M}\right) + \theta v NF - \pi CN,
\end{aligned} \tag{5}$$

where $C(0) = C_0 \geq 0$, $G(0) = G_0 \geq 0$, $F(0) = F_0 \geq 0$, $N(0) = N_0 \geq 0$. Fig. 2 shows the flow chart of the model system (5).

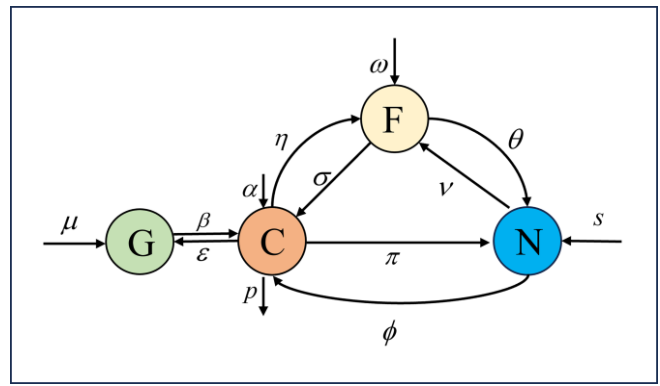

**Fig. 2** Flow chart of the model system (5)

## 3 Model analysis

3.1 Boundedness

The boundedness of the system is given by the following lemma:

**Theorem 3.1.1** If $\alpha + \phi N_{max} + (\beta - \varepsilon)G_{max} > 0$, then solutions of system (5) are bounded in region $\Omega = \{(C, G, F, N) \in R_+^4 : 0 \leq C \leq C_{max}; 0 \leq G \leq G_{max}; 0 \leq F \leq F_{max}; 0 \leq N \leq N_{max}\}$, here $C_{max}$, $G_{max}$, $F_{max}$ and $N_{max}$ are given as $C_{max} = \dfrac{\alpha + \phi N_{max} + (\beta - \varepsilon)G_{max}}{p}$, $G_{max} = \dfrac{\mu}{\varepsilon}$, $F_{max} = \dfrac{K(\omega + \eta \sigma C_{max})}{\omega}$, $N_{max} = M + \dfrac{\theta v M}{s} F_{max}$ and this attracts all solutions initiating from within the positive orthant's interior..

**Proof** According to the comparison theorem we get

$$0 \leq G(t) \leq \frac{\mu}{\varepsilon} = G_{max}(say).$$

From the first equation we get

$$\frac{dC}{dt} \leq \alpha + \phi N + (\beta - u)G_{max} - pC,$$

which gives that

$$\limsup_{t \to \infty} C(t) = \frac{\alpha + \phi N_{max} + (\beta - u)G_{max}}{p} = C_{max}(say).$$

From the third equation we get

$$\frac{dF}{dt} \leq \omega F\left(1 - \frac{F}{K}\right) + \eta \sigma C F + \varepsilon G,$$

which implies that

$$\limsup_{t \to \infty} F(t) = \frac{K(\omega + \eta \sigma C_{max})}{\omega} = F_{max}(say).$$

From the fourth equation we get

$$\frac{dN}{dt} \leq (s + \theta v F)N - \frac{s}{M}N^2,$$

which implies that

$$\limsup_{t \to \infty} N(t) = \frac{M(s + \theta v F_{max})}{s} = M + \frac{\theta v}{s}F_{max} = N_{max}(say).$$

This completes the proof of the boundedness of the system (5).

3.2 Equilibrium points

To solve the equilibrium point of the model, we need to make the right of the model equal to zero to solve all possible equilibrium points. The model system (5) has two non-negative equilibria which

are listed as follows:

(i) $E_1 = (C_1, G_1, 0, 0)$.

(ii) $E_2 = (C_2, G_2, F_2, 0)$.

(iii) $E_3 = (C_3, G_3, 0, N_3)$.

(iv) $E_4 = (C_4, G_4, F_4, N_4)$.

**Theorem 3.2.1** System (5) possesses an equilibrium point $E_1$ provided that the following inequality is satisfied:

$$\alpha + \left(1 - \frac{\beta}{\varepsilon}\right) > 0. \tag{6}$$

**Proof** (i) From equation two we have

$$\frac{dG}{dt} = \mu - \varepsilon G = 0, \tag{7}$$

therefor we have

$$G_1 = G_* = \frac{\mu}{\varepsilon}. \tag{8}$$

When $N = F = 0$, from equation one we have

$$C = \frac{1}{p}\left[\alpha + \left(\frac{\beta}{\varepsilon} - \mu\right)\right]. \tag{9}$$

Thus we get equilibrium point $E_1 = \left(\frac{1}{p}\left[\alpha + \left(\frac{\beta}{\varepsilon} - \mu\right)\right], \frac{\mu}{\varepsilon}, 0, 0\right)$.

**Theorem 3.1.2** System (5) possesses an equilibrium point $E_2$ provided that the following inequality is satisfied:

$$\alpha\varepsilon + \mu(\beta - \varepsilon) > 0. \tag{10}$$

(ii) When and $N = 0$, from equation one and three we have

$$\omega - \frac{\omega}{K}F + \eta\sigma C = 0, \tag{11}$$

$$\alpha + (\beta - \varepsilon)\frac{\mu}{\varepsilon} - \eta CF - pC = 0. \tag{12}$$

Using the value of $F$ from (11) in (12), we get the following quadratic polynomial in $C$

$$\frac{\eta^2 \sigma K}{\omega}C^2 + (\eta K + p)C - \left[\alpha - \left(1 - \frac{\beta}{\varepsilon}\right)\right] = 0. \tag{13}$$

Hence, applying Descartes' rule of signs confirms the existence of a unique positive root

$$C_2 = \frac{\omega\left[-(p+\eta K)+\sqrt{(p+\eta K)^2 + \frac{4\sigma\eta^2 K}{\omega}\left[\alpha+\left(1-\frac{\beta}{\varepsilon}\right)\right]}\right]}{2\sigma\eta^2 K}, \quad (14)$$

if condition $\alpha\varepsilon + \mu(\beta-\varepsilon) > 0$ satisifed. Substitute equation (14) into equation (11), and a unique positive value of $F_2$ can be obtained.

**Theorem 3.2.3** System (5) possesses an equilibrium point $E_3$ provided that the following inequality is satisfied:

$$\varepsilon(\alpha+\phi M) + \mu(\beta-\varepsilon) > 0, \quad (15)$$

$$\varepsilon(sp+\pi\phi M) > \varepsilon(\alpha+\phi M) + \mu(\beta-\varepsilon). \quad (16)$$

(iii) when $F = 0$, from equation one and four of systsem (5) we have

$$\alpha + \phi N + (\beta-\varepsilon)G - pC = 0, \quad (17)$$

$$s - \frac{s}{M}N - \pi C = 0. \quad (18)$$

Using the value of $N$ from (18) in (17) we get

$$C = \frac{s\left[\varepsilon(\alpha+\phi M)+\mu(\beta-\varepsilon)\right]}{\varepsilon(sp+\pi\phi M)}. \quad (19)$$

Using the value of $C$ in (18) we get

$$N = M\left\{1 - \frac{\pi\left[\varepsilon(\alpha+\phi M)+\mu(\beta-\varepsilon)\right]}{\varepsilon(sp+\pi\phi M)}\right\}. \quad (20)$$

**Theorem 3.2.4** System (5) possesses an equilibrium point $E_4$ provided that the following inequality is satisfied:

$$\left[\varepsilon p(sp+\pi\phi M)\right]\left\{\varepsilon\omega p + \eta\sigma\left[\varepsilon\alpha+\mu(\beta-\varepsilon)\right]\right\} > \varepsilon M(p\theta-\eta\sigma\phi)\left[\varepsilon(\alpha+s\pi p)+\mu(\beta-\varepsilon)\right]. \quad (21)$$

(iv) From system (5) we get following equations

$$\alpha + \phi N + (\beta-\varepsilon)\frac{\mu}{\varepsilon} - \eta CF - pC = 0, \quad (22)$$

$$\omega - \frac{\omega}{K}F - \theta N + \eta\sigma C = 0, \quad (23)$$

$$s - \frac{s}{M}N + \theta vF - \pi C = 0. \quad (24)$$

From (22) we have

$$C = \frac{\alpha + \phi N + (\beta - \varepsilon)\frac{\mu}{\varepsilon}}{p + \eta F}. \tag{25}$$

Using (25) in (23) and (24) we get following equations

$$a : \omega - \frac{\omega}{K}F - \theta N + \eta\sigma \frac{\alpha + \phi N + (\beta - \varepsilon)\frac{\mu}{\varepsilon}}{p + \eta F} = 0, \tag{26}$$

$$b : s - \frac{s}{M}N + \theta v F - \frac{\alpha + \phi N + (\beta - \varepsilon)\frac{\mu}{\varepsilon}}{p + \eta F} = 0. \tag{27}$$

To prove the existence of equilibrium points, we now analyze curves $a$ and $b$ separately.

For equation (26):

(i) when $N = 0$, $F = F_2 > 0$ if $\alpha\varepsilon + \mu(\beta - \varepsilon) > 0$.

(ii) when $F = 0$, we have $N = N_a = \dfrac{\varepsilon\omega p + \eta\sigma[\varepsilon\alpha + \mu(\beta - \varepsilon)]}{\varepsilon(p\theta - \eta\sigma\phi)} > 0$ if

$$p\theta - \eta\sigma\phi > 0. \tag{28}$$

(iii) By differentiating equation (26) with respect to $F$, we obtain:

$$\frac{dN}{dF} = \frac{\varepsilon\omega(p + \eta F)^2 + \sigma\eta^2[\varepsilon(\alpha + \phi N) + \mu(\beta - \varepsilon)]K}{\varepsilon K(p + \eta F)[\eta\sigma\phi - \theta(p + \eta F)]} < 0. \tag{29}$$

For equation (27):

(i) when $F = 0$, we have $N = N_b = \dfrac{M[\varepsilon(\alpha + s\pi p) + \mu(\beta - \varepsilon)]}{\varepsilon p(sp + \pi\phi M)} > 0$ if

$$\varepsilon(\alpha + s\pi p) + \mu(\beta - \varepsilon) > 0. \tag{30}$$

(ii) when $N = 0$, we get following equation in $F$

$$v\eta\theta F^2 + s\varepsilon(\eta + vp)F + [\varepsilon(sp - \alpha) - (\beta - \varepsilon)] = 0, \tag{31}$$

applying Descartes' rule of signs confirms the existence of a unique negative root $F = F_b < 0$ if

$$\varepsilon(sp - \alpha) - (\beta - \varepsilon) < 0 \tag{32}$$

(iii) Calculate the derivative of (27) with respect to $F$ we get

$$\frac{dN}{dF} = \frac{\nu\varepsilon\theta(p+\eta F)^2 + \eta[\varepsilon(\alpha+\phi N)+\mu(\beta-\varepsilon)]}{\varepsilon(p+\eta F)[\phi M + s(p+\eta F)]} M > 0. \tag{33}$$

This indicates that there is a unique intersection between the two curvesn equilibrium point $E_4$, at this time $N_a > N_b$ must be satisfied, i.e

$$[\varepsilon p(sp+\pi\phi M)]\{\varepsilon\omega p + \eta\sigma[\varepsilon\alpha + \mu(\beta-\varepsilon)]\} > \varepsilon M(p\theta - \eta\sigma\phi)[\varepsilon(\alpha+s\pi p)+\mu(\beta-\varepsilon)].$$

## 4 Stability analysis

We discuss the stability of equilibria $E_1$, $E_2$ and $E_3$ by finding the sign of the eigenvalues of Jacobian matrix corresponding to each equilibrium. The Jacobian matrix for model system (5) is given as follows:

$$J = \begin{bmatrix} -\eta F - p & \beta - \varepsilon & -\eta C & \phi \\ 0 & -\varepsilon & 0 & 0 \\ \eta\sigma F & 0 & \omega - \frac{2\omega}{K}F - \theta N + \eta\sigma C & -\theta F \\ -\pi N & 0 & \nu\theta N & s - \frac{2s}{M}N + \nu\theta F - \pi C \end{bmatrix}. \tag{34}$$

Define $J_i (i=1,2,3,4)$ as the equilibrium Jacobian.

**Theorem 4.1** (i) $E_1$ is inherently unstable under all conditions.

(ii) $E_2$ is always stable in $G$ direction and locally stable (unstable) manifold in $N$ direction provideds $s + \nu\theta F_2 - \pi C_2$ is negative (positive). Asymptotically stable in $C-F$ directions when

$$\omega - \frac{2\omega}{K}F_2 + \eta\sigma C_2 < \min\left\{\eta F_2 + p, \frac{\sigma\eta^2 F_2 C_2}{\eta F_2 + p}\right\}.$$

(iii) $E_3$ is always stable in $G$ direction, whereas $E_3$ is locally stable (unstable) manifold in $F$ direction provideds $\omega - \frac{2\omega}{K}F_2 + \eta\sigma C_2$ is negative (positive). asymptotically stable in n $C-F$ directions when $s - \frac{2s}{M} - \pi C_3 < \min\left\{p, \frac{\pi\phi N_3}{p}\right\}$.

**Proof** (i) The eigenvalues of the Jacobian matrix $J_1$ are $-p$, $-\varepsilon$, $\omega + \eta\sigma\left[\frac{1}{p} + \left(\frac{\beta}{\varepsilon} - \mu\right)\right]$ and $-\pi\left[\frac{1}{p} + \left(\frac{\beta}{\varepsilon} - \mu\right)\right]$ seperately. $\omega + \eta\sigma\left[\frac{1}{p} + \left(\frac{\beta}{\varepsilon} - \mu\right)\right] > 0$ whenever $E_1$ exist.

(ii) The eigenvalues of the Jacobian matrix $J_2$ in $G$ and $N$ directions are $-\varepsilon$ and $s+v\theta F_2 - \pi C_2$, therefore $E_2$ is always stable in $G$ direction, whereas locally stable (unstable) manifold in $N$ direction provideds $s+v\theta F_2 - \pi C_2$ is negative (positive). The other two eigenvalues are solutions to unary quadratic equation $y^2 - \left(\omega + \eta\sigma C_2 - \eta F_2 - p - \frac{2\omega}{K}F_2\right) + \left[\sigma\eta^2 F_2 C_2 - (\eta F_2 + p)\left(\omega - \frac{2\omega}{K}F_2 + \eta\sigma C_2\right)\right] = 0$. When the trace is less than zero and the value of the determinant is greater than zero, there are negative eigenroots or negative real part, ie

$$\omega - \frac{2\omega}{K}F_2 + \eta\sigma C_2 < \min\left\{\eta F_2 + p, \frac{\sigma\eta^2 F_2 C_2}{\eta F_2 + p}\right\}. \tag{35}$$

Thus the two eigenroots are asymptotically stable when (30) satisfied.

(iii) The eigenvalues of the Jacobian matrix $J_3$ in $G$ and $F$ directions are $-\varepsilon$ and $\omega - \theta N_2 + \eta\sigma C_2$, therefore $E_3$ is always stable in $G$ direction whereas $E_3$ is locally stable (unstable) manifold in $F$ direction provideds $\omega - \frac{2\omega}{K}F_2 + \eta\sigma C_2$ is negative (positive). The other two eigenvalues are solutions to unary quadratic equation $y^2 - \left(s - p - \frac{2s}{M} - \pi C_3\right) + \left[\pi\phi N_3 - p\left(s - \frac{2s}{M} - \pi C_3\right)\right] = 0$. When the trace is less than zero and the value of the determinant is greater than zero, there are negative eigenroots or negative real part, ie

$$s - \frac{2s}{M} - \pi C_3 < \min\left\{p, \frac{\pi\phi N_3}{p}\right\}. \tag{36}$$

Thus the two eigenroots are asymptotically stable when (31) satisfied.

**Theorm 4.2** Local asymptotic stability of the interior equilibrium point $E_4$ is guaranteed if the subsequent condition holds:

$$A_1 A_2 - A_3 > 0. \tag{37}$$

(iv) **Proof** Evaluating the Jacobian matrix at $E_4$ we derive:

$$J_4 = \begin{bmatrix} -\eta F_4 - p & \beta - \varepsilon & -\eta C_4 & \phi \\ 0 & -\varepsilon & 0 & 0 \\ \eta\sigma F_4 & 0 & \omega - \frac{2\omega}{K}F_4 - \theta N_4 + \eta\sigma C_4 & -\theta F_4 \\ -\pi N_4 & 0 & v\theta N_4 & s - \frac{2s}{M}N_4 + v\theta F_4 - \pi C_4 \end{bmatrix}.$$

The characteristic equation of $J_4$ is

$$(\psi + \varepsilon)(\psi^3 + A_1\psi^2 + A_2\psi + A_3) = 0,$$

where,

$$A_1 = \eta F_4 + p + \frac{2\omega}{K}F_4 + 2\nu\theta F_4 + \frac{s}{M}N_4,$$

$$A_2 = \frac{\omega}{K}F_4 \cdot \frac{s}{M}N_4 + (p + \eta F_4)\left(\frac{\omega}{K}F_4 + \frac{s}{M}N_4\right),$$

$$A_3 = \frac{s}{M}\eta\theta(\pi C_4 + \nu\sigma\phi)N_4^2 F_4 + \left[\frac{\omega}{K}\pi\phi + \nu\theta^2(\eta F_4 + p)\right]N_4 F_4 + \frac{\omega}{K}\cdot\frac{s}{M}\eta^2\sigma N_4 F_4^2.$$

Here, it is straightforward to see that $A_1$, $A_2$ and $A_3$ are positive. The Routh–Hurwitz criterion therefore simplifies to $A_1 A_2 - A_3 > 0$.

**Theorm 4.3** Global stability of the interior equilibrium $E_4$ inside the region of attraction is guaranteed if the subsequent conditions hold:

$$\max\left\{\frac{\eta K}{m_1\omega}(m_1\sigma - C_{\max})^2, \frac{(\beta - \varepsilon)^2}{\varepsilon}\right\} < 2(p + \eta F_4). \tag{38}$$

**Proof** To establish the global stability of interior equilibrium $E_4$, we employ Lyapunov's method by selecting a positive definite function as:

$$V = \frac{1}{2}(C - C_4)^2 + \frac{1}{2}(G - G_4)^2 + m_1\left(F - F_4 - F_4 \ln\frac{F}{F_4}\right) + m_2\left(N - N_4 - N_4 \ln\frac{N}{N_4}\right),$$

where $m_1$ and $m_2$ represent positive constants subject to appropriate calibration.

Differentiating $V$ with respect to $t$ along the solution path of model system (5) yields:

$$\frac{dV}{dt} = -(p + \eta F_4)(C - C_4)^2 - \varepsilon(G - G_4)^2 - \frac{\omega}{K}(F - F_4)^2 - \frac{s}{M}m_2(N - N_4)^2$$
$$+ (\phi - \pi m_2)(N - N_4)(C - C_4) + (\beta - \varepsilon)(C - C_4)(G - G_4)$$
$$+ \eta(m_1\sigma - C)(C - C_4)(F - F_4) + \theta(\nu m_2 - m_1)(F - F_4)(N - N_4),$$

choosing $m_2 = \dfrac{\phi}{\pi}$ and $m_1 = \nu m_2 = \dfrac{\nu\phi}{\pi}$ we get

$$\frac{dV}{dt} = \left[-\frac{p + \eta F_4}{2}(C - C_4)^2 + \eta(m_1\sigma - C)(F - F_4)(C - C_4) - \frac{m_1\omega}{K}(F - F_4)^2\right]$$
$$+ \left[-\varepsilon(G - G_4)^2 + (\beta - \varepsilon)(C - C_4)(G - G_4) - \frac{p + \eta F_4}{2}(C - C_4)^2\right] - \frac{m_2 s}{K}(N - N_4)^2$$

$$\leq \left[-\frac{p+\eta F_4}{2}(C-C_4)^2 + \eta(m_1\sigma - C_{max})(F-F_4)(C-C_*) - \frac{m_1\omega}{K}(F-F_4)^2\right]$$

$$+\left[-\varepsilon(G-G_4)^2 + (\beta-\varepsilon)(C-C_4)(G-G_4) - \frac{p+\eta F_4}{2}(C-C_4)^2\right] - \frac{m_2 s}{K}(N-N_4)^2.$$

$\frac{dV}{dt}$ is negative when $\begin{cases} \frac{\eta^2 K}{m_1\omega}(m_1\sigma - C_{max})^2 < 2(p+\eta F_4) \\ \frac{(\beta-\varepsilon)^2}{\varepsilon} < 2(p+\eta F_4) \end{cases}$ statisted and the intersection of these two conditions is:

$$\max\left\{\frac{\eta K}{m_1\omega}(m_1\sigma - C_{max})^2, \frac{(\beta-\varepsilon)^2}{\varepsilon}\right\} < 2(p+\eta F_4).$$

Now, we observe that $\frac{dV}{dt}$ satisfies negative definiteness inside the attraction region '$X$', provided that condition (38) is met.

## 5 Parameter estimation

For parameter estimation, China's $CO_2$ emission data is utilized, encompassing emissions from fossil fuel combustion and land-use change while excluding other carbon emission sources (Ritchie and Roser. 2020), GDP (Kylc. 2024), human population (NBSC. 2024) and forest area (Kylc. 2023) from 2000 to 2022. The natural rate of growth of atmospheric $CO_2$ we take 1.68 ppm per year (Verma and Verma. 2021). For the 2000-2022 period, the average annual per capita anthropogenic $CO_2$ emission rate is 1.46 billion tons, equivalent to 0.08 parts per million (ppm) per million people annually. As same as (Caetano. 2008) we take $\varepsilon = 0.0008$. From 2001 to 2010, the atmospheric lifetime of carbon dioxide typically ranges from 30 to 95 years (Verma and Misra. 2018), accordingly, the natural sink rate of atmospheric $CO_2$ is 0.016. The average growth rate of GDP $\mu$ in China from 2000 to 2022 is 0.02145. According to (Kylc. 2023), the growth rate of forest area $\omega$ we take 0.06133. The average intrinsic growth rate of population is approximately 0.00529 during this period of time. We take the mortality rate of the population caused by global warming as 0.00005 per ppm per year (Devi and Gupta. 2019). Because there is no actual data to support it, the rate of deforestation and the depelation rate of carbon dioxide due to forest we take 0.0004 and 0.0000001 respectively. Forests absorb carbon dioxide through photosynthesis to promote themselves at rate 0.01 per million heacter per year (Misra and Verma. 2013; Lata and Misra; 2015 ). We assume that the rate of human population increase due to forest is proportional to the consumption of forest resources by a rate 0.001 (Aristide and Mayengo. 2014). As same as (Caetano et al. 2011 ) we take $\beta = 0.0003$. Based on the above assumptions and conjectures, the following parameter values will be used in the subsequent numerical simulation.

$$\alpha = 1.68, \phi = 0.008, \beta = 0.0003, \eta = 0.0000001,$$
$$p = 0.016, \mu = 0.02145, \varepsilon = 0.0008, \omega = 0.06133,$$
$$s = 0.00529, M = 1720, K = 11000, \theta = 0.0004,$$
$$\sigma = 0.01, \nu = 0.001, \pi = 0.00005.$$

## 6 Numerical simulation

For the confirmation and graphical representation of analytical findings, we simulated model system (5) with the parameter values listed in Section 5. The numerical simulation was carried out utilizing MATLAB R2023a. The interior equilibrium components are obtained as: $E_4(130.9959, 26.8125, 3607.3559, 59.5748)$. The Jacobian matrix eigenvalues at equilibrium $E_4$ are calculated as $-0.0084$, $-0.0102$, $-0.0262$ and $-0.0008$—all negative—thus establishing the local asymptotic stability of equilibrium point $E_4$. For the data mentioned above, the solution trajectories of the model system (5) have been plotted in Fig. 3 and Fig. 4 with different initial conditions. As observed, all trajectories starting inside the region of attraction tend toward equilibrium point $E_4$, demonstrating the nonlinear stability of interior equilibrium in the $C-G-F$ and $N-F-C$ spaces.

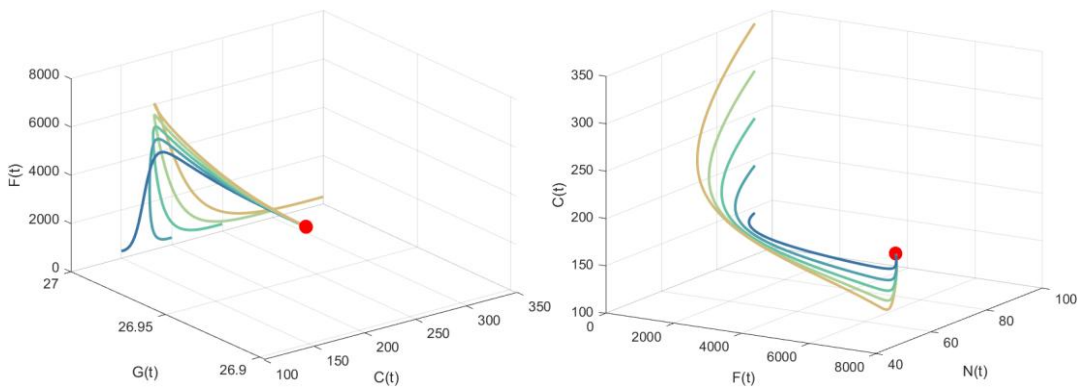

**Fig. 3** Global Stability inside the *C-G-F* space  **Fig. 4** Global Stability inside the *N-F-C* space

Fig.5 and Fig.6 depict highly significant results for the investigated dynamical system. These graphs are plotted to observe the temporal changes in the concentration of carbon dioxide $C(t)$, human population $N(t)$ and forest area $F(t)$. An examination of the variations is performed with regard to various parameter values $\phi$ and $\pi$. By comparison, the values of other parameters stay unchanged, as demonstrated in Table 1. From Fig. 5, we observe that as the anthropogenic carbon dioxide emission rate ($\phi$) increases from 0.005 to 0.006, the carbon dioxide concentration rises from 122.1270 ppm to 128.1255 ppm, while the human population declines from 129.2128 million to 126.8650 million. When $\phi$ increases further from 0.006 to 0.007, the $CO_2$ concentration increases to 134.0093 ppm, and the human population decreases to 124.5806 million. Consistent trends in $CO_2$ concentration and human population persist with additional increases in the parameter $\phi$.

From Fig. 6, as the parameter $\pi$ increases from 0.0005 to 0.0006, the carbon dioxide concentration drops from 139.7815 ppm to 137.1956 ppm, while the forest area expands from 1707.6322 million hectares to 1980.1586 million hectares. When $\pi$ further increases from 0.0006 to 0.0007, the $CO_2$ concentration decreases to 134.7915 ppm, and the forest area grows to 2264.8575 million hectares. Consistent patterns in $CO_2$ concentration and forest area persist with additional increments in parameter $\pi$. Thus, according to these numerical findings, it can be deduced that the rate of anthropogenic carbon dioxide emissions elevates atmospheric carbon dioxide concentration and exerts a certain level of detriment to human health. Conversely, forests are capable of absorbing carbon dioxide via photosynthesis, a process that not only fosters their own growth but also markedly decreases the atmospheric carbon dioxide concentration.

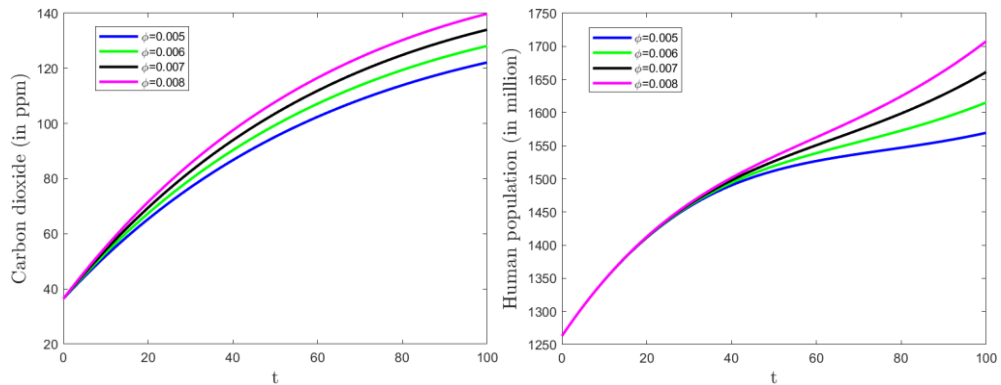

**Fig. 5** Time series graph of $C$ and $N$ for different value of parameter $\phi$

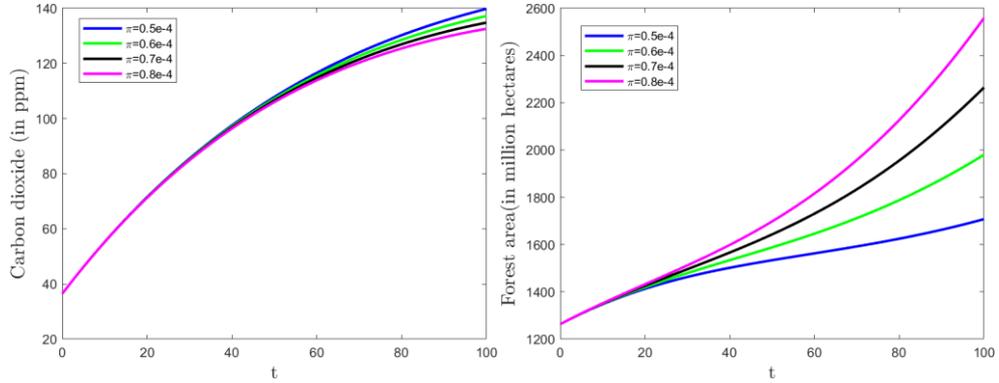

**Fig. 6** Time series graph of $F$ and $N$ for different value of parameter $\pi$

## 7 Sensitivity analysis

To elucidate the impact of model parameters on system dynamics, we implemented a global sensitivity analysis (GSA) via the partial ranking correlation coefficient (PRCC) approach, which leverages Latin Hypercube sampling Monte Carlo simulation (LHS). As documented by (Bidah et al. 2020), this methodology facilitates the evaluation of individual parameter fluctuations on the aggregate model response. A positive PRCC value denotes a direct dependency between model parameters and their outputs, such that an increment in parameter values typically elicits a pronounced rise in model output, whereas a decrement generally results in output reduction. Conversely, (Fanuel et al. 2023) illustrate that a negative PRCC value implies an inverse relationship: increasing parameter values correspond to decreasing model outputs, and conversely, decreasing parameters yield increasing outputs. The baseline values and interval of the parameters for sensitivity analysis are provided in Table 1.

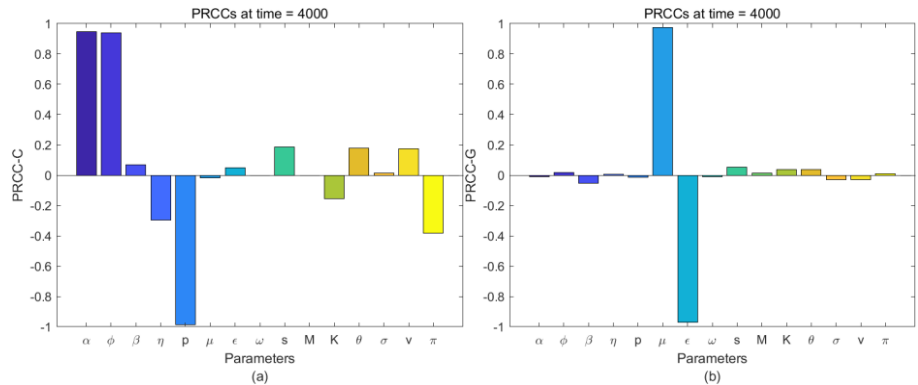

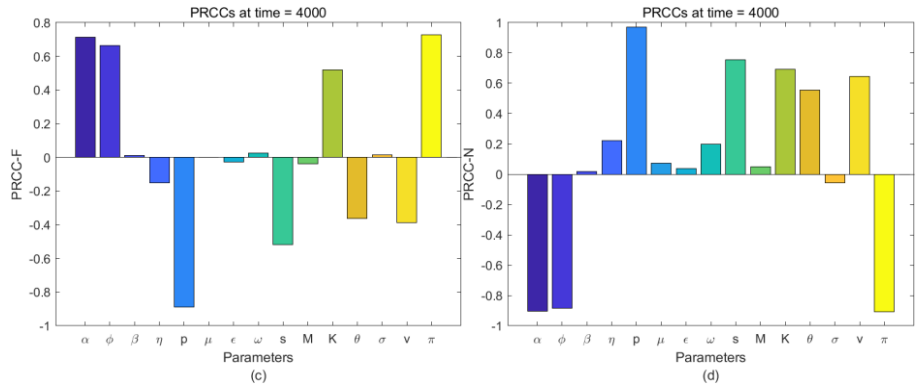

**Fig .9** The PRCC value of each parameter with respect to each compartment at t=4000 respectively

From Fig. 9(a), parameters exerting positive effects on the compartment $C(t)$ are identified as $\alpha$ and $\phi$, while $p$ demonstrates the most pronounced negative impact on the compartment $C(t)$. Fig. 9(b) reveals that $\mu$ and $\varepsilon$ respectively exert the strongest positive and negative effects on compartment $G(t)$, with other parameters showing negligible influence. In Fig. 9(c), parameters contributing positively to compartment $F(t)$ include $K$, $\pi$, $\alpha$ and $\phi$, whereas $s$ and $p$ exhibit the most significant negative effects; other parameters have minimal impact. Fig. 9(d) shows that $p$, $s$, $K$, $v$ and $\theta$ positively affect compartment $N(t)$, while $\pi$, $\alpha$ and $\phi$ exert negative effects, with remaining parameters having insignificant influence.

**Table 1**. Parameter with their basline.

| Paramater | Description | Basline | Interval |
|---|---|---|---|
| $\alpha$ | Natural $CO_2$ emission rate | 1.68 | [1.521 1.848] |
| $\phi$ | Anthropogenic $CO_2$ emission rate | 0.008 | [0.0072 0.0088] |
| $\varepsilon$ | GDP-driven $CO_2$ decay rate | 0.0008 | [0.00072 0.00088] |
| $p$ | $CO_2$ decay rate | 0.016 | [0.0144 0.0176] |
| $\eta$ | $CO_2$ depletion coefficient due to forest area | 0.000001 | [0.0000009 0.0000011] |
| $\mu$ | GDP growth rate | 0.02145 | [0.019305 0.023595] |
| $\omega$ | Inherent growth rate of forest area | 0.06133 | [0.055197 0.067463] |
| $K$ | Forest area carrying capacity | 11000 | [10000 12000] |
| $\theta$ | Rate of forest loss | 0.0004 | [0.00036 0.00044] |
| $\sigma$ | Growth rate of forest area for $CO_2$ absorption | 0.01 | [0.009 0.011] |
| $\beta$ | $CO_2$ emissions in the Process of GDP Growth | 0.0003 | [0.00027 0.00033] |
| $s$ | Inherent growth rate of human population | 0.00529 | [0.004761 0.005819] |
| $M$ | Human population carrying capacity | 1720 | [1542 1892] |
| $\nu$ | Forest area-related human population growth rate | 0.001 | [0.0009 0.0011] |
| $\pi$ | Mortality rate coefficient from elevated $CO_2$ | 0.00005 | [0.000045 0.000055] |

## 8 Optimal control

The rise in carbon emissions is intrinsically linked to China's rapid GDP growth. However, once economic development reaches a certain threshold, these economic resources can be effectively harnessed to implement technological interventions and concerted efforts aimed at absorbing and mitigating atmospheric carbon dioxide concentrations, we denote these measures by $u$. Nonetheless, a significant budgetary allocation is needed to fund the expenses related to these measures. Therefore, in terms of project implementation, it is necessary to formulate a cost-optimal intervention strategy, with an implementation speed that is sufficient to carry out adequate measures and actions while also minimizing the implementation cost. Given that $u$ is not a constant but a Lebesgue measurable function over the finite time interval $[0, t_f]$, the model system can be rewritten as:

$$\begin{aligned}
\frac{dC}{dt} &= \alpha + \phi N + (\beta - u(t))G - \eta CF - pC, \\
\frac{dG}{dt} &= \mu - u(t)G, \\
\frac{dF}{dt} &= \omega F\left(1 - \frac{F}{K}\right) - \theta NF + \eta \sigma CF, \\
\frac{dN}{dt} &= sN\left(1 - \frac{N}{M}\right) + \theta \nu NF - \pi CN,
\end{aligned} \tag{39}$$

where $C(0) = C_0 \geq 0$, $G(0) = G_0 \geq 0$, $F(0) = F_0 \geq 0$, $N(0) = N_0 \geq 0$.

To minimize the objective cost function, we employ Pontryagin's maximum principle (Pontryagin et al., 1962). The specific form of the objective cost function is detailed as follows:

$$J = \min_{u} \int_0^{t_f} \left[AC(t) + \frac{B}{2}u^2(t)\right]dt, \tag{40}$$

Here, the parameters $A$ and $B$ represent the weighting parameters of the function (40). The term $\frac{B}{2}u^2(t)$ characterizes the cost incurred by the measures and actions. Bubject to model (39), we seek optimal control $u_*(t)$ such that

$$J(u_*(t)) = \min_{u(t)\in\Theta} J(u(t)), \tag{41}$$

where the control set is denoted by $\Theta = \{u(t): 0 \leq u(t) \leq u_{\max} \text{ for } t \in [0, t_f]\}$.

**Theorm 8.1** On a fixed interval $[0, t_f]$, an optimal control $u_* \in \Theta$ exists to minimize the objective function (40) under the constraint of system (39).

Referring to Aldila (2020), the optimal control problem under our consideration should comply with the conditions as follows:

1. Let $x = (C(t), G(t), F(t), N(t))$, for a given initial value $x_0$, the set $\{(x_0, u)\}$ composed of control variable $u$ and the solutions of the state equation that satisfy the initial conditions are non-empty.

2. $\Theta$ should be closed and convex. System (39) is a function of the control variable $u$, and the coefficients of the objective function depend on time and state variables.

3. $D = AC(t) + \frac{B}{2}u^2(t)$ is convex on $\Theta$ and satisfies $D \geq f(u)$, here $f(u)$ is a continuous function and satisfies condition $\lim_{|u|\to\infty} \frac{f(u)}{|u|} = \infty$. Note that $|\cdot|$ represent the norm.

**Proof** According to the proof of the boundedness of the model, we know that $x$ is bounded, that is $x_{\max} = (C_{\max}, G_{\max}, F_{\max}, N_{\max})$. As long as $u$ is bounded within $\Theta$, the solutions of system (39) are always bounded. Therefore, condition one is satisfied.

According to the definition, given a control set $\Theta$, where $u \in [0, 1]$, thus set $\Theta$ is closed. According to the definition of a convex set, let $\Theta$ be a set. If for any $x_1, x_2 \in \Theta$ and any real number $\delta \in [0, 1]$, we have that $\delta x_1 + (1-\delta)x_2 \in \Theta$, then $\Theta$ is called a convex set. Therefor, $\delta x_1 + (1-\delta)x_2 \in \Theta$ implying $\Theta$ is convex. System (39) can obviously be expressed as a function of the control variable $u$. The coefficients of the objective function $A$ and $B$ depend on time $t$ and state variables $C(t)$. Hence, we fulfilled 2.

The integrand $D = AC(t) + \frac{B}{2}u^2(t)$ is convex due to the quadratic form of $u$. Further more, $D = AC(t) + \frac{B}{2}u^2(t) \geq \frac{B}{2}u^2(t) = f(u)$. Obviously, $f(u)$ is continuous and satisfies condition $\lim_{|u|\to\infty} \frac{f(u)}{|u|} = \infty$. Hence, we fulfilled 3. Therefor, an optimal control $u_* \in \Theta$ exists to minimize the objective function (40) under the constraint of system (39) over the fixed interval $[0, t_f]$.

Employing Pontryagin's maximum principle to characterize the optimal control, the Hamiltonian is given in the following way.

$$H(C, G, F, N, u, u_1, u_2, u_3, u_4) = AC(t) + \frac{B}{2}u^2(t)$$
$$+ u_1 \left[\alpha + \phi N + (\beta - u)G - \eta CF - pC\right]$$
$$+ u_2 (\mu - uG)$$
$$+ u_3 \left[\omega F\left(1 - \frac{F}{K}\right) - \theta NF + \eta \sigma CF\right]$$
$$+ u_4 \left[sN\left(1 - \frac{N}{M}\right) + \theta \nu NF - \pi CN\right],$$

where $u_i$ $(i = 1, 2, 3, 4)$ are the adjoint variable to be determined by solving the following equations

$$u_1' = -\frac{\partial H}{\partial C} = -A + u_1(\eta F + p) - u_3 \eta \sigma F + u_4 \pi N,$$
$$u_2' = -\frac{\partial H}{\partial G} = u_1(u - \beta) + u_2 u,$$
$$u_3' = -\frac{\partial H}{\partial F} = u_1 \eta C - u_3 \left[\omega\left(1 - \frac{2F}{K}\right) - \theta N + \eta \sigma C\right] - u_4 \theta \nu N,$$
$$u_4' = -\frac{\partial H}{\partial N} = -u_1 \phi + u_3 \theta F - u_4 \left[s\left(1 - \frac{2N}{M}\right) + \theta \nu F - \pi C\right],$$

along with transversality conditions

$$u_1(t_f) = u_2(t_f) = u_3(t_f) = u_4(t_f) = 0, \tag{42}$$

and

$$u_* = \max\left\{0, \min\left\{u_{\max}, \frac{(u_1 + u_2)G}{B}\right\}\right\}. \tag{43}$$

## 9 Numerical simulation of optimal control

To demonstrate the optimal mitigation strategies for the control of future CO₂ level, the optimality system (39) is solved numerically by choosing the upper limit of the control $u_{\max} = 0.008$, weight parameters $A = 0.0001$ and $B = 10$, final time $t_f = 100$. We utilized the forward-backward sweep method to numerically solve the optimality system corresponding to the parameter values in Section 5. First, we initialized the control variable with reasonable guesses. The state equations were then integrated forward in time using the fourth-order Runge-Kutta method, while the adjoint equations were solved backward in time. The control was updated iteratively until convergence. This process was repeated until the attainment of the desired convergence. Initial states are set as $C_0 = 130$, $G_0 = 0.121$, $F_0 = 1003$, $N_0 = 80$. The solution trajectories for the concentration of carbon dioxide $C(t)$, the forest area $F(t)$, and the human population $N(t)$, both under the dynamic optimal control and in the absence of control strategies, are depicted in Fig. 10. This figure clearly shows the significant reduction in carbon dioxide under the time-dependent optimal control. It is plainly evident from these graphs that the optimal control strategy outperforms the strategy without control, effectively demonstrating its superiority. The solution trajectories for atmospheric $CO_2$ concentration $C(t)$, GDP $G(t)$, forest area $F(t)$, and human population $N(t)$—under both dynamic optimal control and uncontrolled conditions—are shown in Fig. 10. The figure exhibits a marked decrease in $CO_2$ concentration under time-varying optimal control. Notably, the plots clearly demonstrate that the optimal control strategy outperforms the uncontrolled scenario, effectively confirming its superiority.

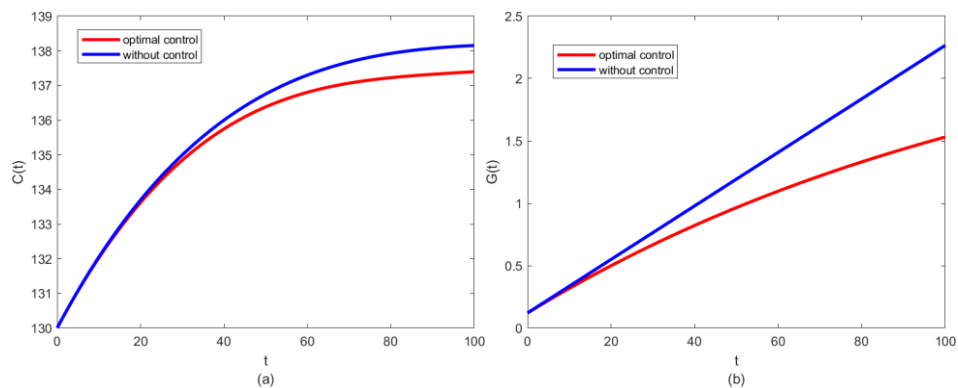

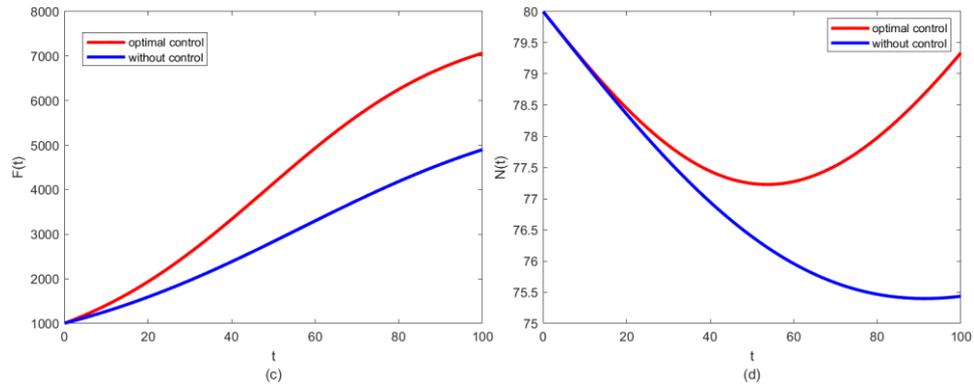

**Fig .10** Graph trajectories with and without optimal control for: (a) carbon dioxide, (b) GDP, (c) forest area, (d) human population

## 10 Conclusion

With the rapid advancement of human society, an increasing amount of greenhouse gases, predominantly carbon dioxide—are being emitted into the Earth's atmosphere, giving rise to severe global warming and climate change. For the sustainable development of human society, understanding the dynamic behavior of atmospheric carbon dioxide is crucial for mitigating these environmental challenges. China currently contributes around one-third of global carbon emissions, a scale closely associated with its rapid GDP growth. As such, the study of China's carbon emission trends is of great significance for global climate change mitigation efforts. In this research, we put forward and analyze a nonlinear mathematical model that establishes a correlation between $CO_2$ emissions and GDP, forest area, as well as population size.

The proposed model is a valuable tool for predicting the long-term impact of China's GDP on the evolution of atmospheric carbon dioxide. The boundedness of the system is verified using the comparison theorem. The conditions for the asymptotic stability of the four equilibrium points are obtained by the eigenvalues of the Jacobian matrix and the Hurwitz criterion. The condition for the global stability of the coexistence equilibrium point is obtained by constructing the Lyapunov function. To validate the model, we simulate the model parameter values in Table 1. Relevant numerical results show that although the increase in GDP is accompanied by a rise in carbon emissions, some measures can be taken by utilizing the growing GDP to mitigate carbon dioxide emissions.

The sensitivity of the compartments with respect to the parameters is analyzed by means of the PRCC (Partial Rank Correlation Coefficient) and the Latin hypercube sampling test. The results indicate that the parameter $\beta$ has a positive impact on carbon dioxide emissions, that is to say, the growth of GDP will be accompanied by an increase in carbon dioxide emissions. The parameter $\varepsilon$ has a negative impact on carbon dioxide, this demonstrates that monetary approaches serve as a rather effective tool for carbon dioxide governance.

Various human activities have led to the release of a large amount of greenhouse gases into the atmosphere, resulting in the rise of atmospheric temperature and climate change. High speed GDP growth is often accompanied by continuous industrialization and the burning of fossil fuels. The carbon emissions from the burning of fossil fuels and industrialization account for 90% of the global total carbon emissions. China is the largest developing country in the world and the country with the largest carbon emissions, accounting for one-third of the global total carbon emissions. Therefore, China has an unshirkable responsibility in dealing with global warming and climate change. Analysis based on models shows that there is a close relationship between China's GDP and carbon emissions. Through reasonable strategies, a win win situation of economic growth and carbon reduction can be achieved. To ensure stable GDP growth while reducing carbon emissions requires the joint efforts of the government, enterprises, and all sectors of society. By means of technological innovation, industrial structure adjustment, energy structure optimization, policy guidance, and the improvement of public awareness, China can be promoted to transform into a low-carbon economy and achieve sustainable development.

**Acknowledgments** This work was supported by the Fundamental Research Funds for the Central Universities (31920250001; 31920250031), the Gansu Provincial Education Department's Graduate Student 'Innovation Star' Project (2025CXZX-249), and the Leading Talents Project of State Ethnic Affairs Commission of China and the Innovation Team of Ecosystem Restoration Modeling Theory and Application of Northwest Minzu University (10017632).